\documentclass{amsart}

\usepackage[all]{xy}
\usepackage{graphicx}
\usepackage{psfrag}
\usepackage{amsthm}
\usepackage{amscd}
\usepackage{amsfonts}
\usepackage{amstext}
\usepackage{amssymb}
\usepackage{amsmath}
\usepackage{amsxtra}
\usepackage{plain}
\usepackage[all]{xy}

\newcommand{\comm}[1]{}

\numberwithin{equation}{section}
\let\newpf\proof  
\newenvironment{pf}{\newpf\proofname}{\qed\endtrivlist}

\let\newpf\proof

\def\be{\begin{equation}}
\def\ee{\end{equation}}

\def\ba{\begin{align}}
\def\ea{\end{align}}

\newtheorem{prop}{{\bf Proposition}}

\newtheorem{lem}[prop]{\bf Lemma}
\newtheorem{thm}[prop]{{\bf Theorem}}

\newtheorem{rem}{\bf Remark}

\newcommand{\N}{\mathbb{N} }

\newcommand{\R}{\mathbb{R} }

\newcommand{\ind}{\rm{ind}}

\newcommand{\End}{\rm{End}}

\newcommand{\id}{\rm{id}}

\begin{document}
\title{Morse inequalities for manifolds with boundary}
\author{M.E. ZADEH}
\address{Mathematisches Institut
Georg-August-Universität
Göttingen, Germany\and Institute for Advanced Studies in Basic Sciences (IASBS)
Zanjan, Iran}
\email{zadeh@uni-math.gwdg.de}
\maketitle

\begin{abstract}
The aim of this paper is to provide a proof for a version of 
Morse inequality for manifolds with boundary. Our main results are certainly 
known to the experts on Morse theory, nevertheless 
it seems  necessary to write down a complete proof for it. Our proof is analytic and 
is based on J. Roe's account of Witten's approach to Morse Theory. 
\end{abstract}
\section{Introduction}
Let $M$ be a smooth manifold with boundary $\partial M\neq\emptyset$ and 
let $f\colon M\to R$ be a smooth function whose  critical points are 
either isolated nondegenerate located in $M\verb+\+\partial M$ or  
located on the boundary such that each connected component of the boundary 
is a Bott manifold of index $0$ or $1$ as we shall describe.  
Let $\partial M=N_+\sqcup N_-$ be a disjoint union of closed manifolds 
such that $N_+$ and $N_-$ are Bott 
submanifolds with indices, respectively, $0$ and $1$. So in a collar 
neighborhood $[0,1)\times N_+$ we have $f(u,y)=\frac{1}{2}u^2$ while in a collar 
neighborhood $[0,1)\times N_-$ we have $f(u,y)=-\frac{1}{2}u^2$. 
In both case $f$ satisfies the von Neumann condition 
\begin{equation}\label{vonboun}
\frac{\partial f}{\partial u}\,(0,y)=0.
\end{equation}
Let $N_+=N_{a+}\sqcup N_{r+}$ and $N_-=N_{a-}\sqcup N_{r-}$ be disjoint union 
of closed manifolds. The subscripts "$a$" and "$r$" 
refer respectively to absolute and relative 
boundary condition that we explain in below. In sequel 
we consider a riemannian metric on $M$ which is Euclidian around non 
degenerate critical points with respect to coordinates provided by Morse 
lemma. Moreover in the collar neighborhood $(1,0]\times\partial M$ it 
is assumed to take the product form $g=d^2u+g_0$ where $g_0$ is a 
riemannian metric on $\partial M$. 
So the flow 
of $\nabla f$, the gradient of $f$ with respect to this metric, is defined 
for $t\in \R$ and limit points of 
each integral curve is either a non degenerate critical point in the interior 
of $M$ or a  point on $\partial M_i$. Put  
\begin{align*}
&\beta_k=\dim H^k_{dr}(M,N_r),\\
&\gamma_k=\dim H^k_{dr}(N_{r-}),\\
&\eta_k=\dim H^k_{dr}(N_{a+}).
\end{align*}
We denote by $c_k$ the number 
of nondegenerate critical points of Morse index $k$. The aim of this note is 
to provide a proof for the following theorem by using the Witten 
approach to Morse theory \cite{Wi1}
\begin{thm}\label{morsein}
The following inequalities hold for $0\leq k\leq n$ 
\begin{equation*}
\mu_k-\mu_{k-1}+\dots+(-1)^k\mu_0\geq\beta_k-\beta_{k-1}+
\dots+(-1)^k\beta_0,
\end{equation*}
where 
\[\mu_k=c_k+\eta_k+\gamma_{k-1}.\] 
The equality holds for $k=n$.
\end{thm}
Notice that the equality for $k=n$ reads 
\begin{equation}\label{eure}
\sum_{k=0}^n(-1)^{n-k}c_k=\chi(M,N_r)+\chi(N_{a+})-\chi(N_{r-})
\end{equation}
As a special case let $\dim M$ be an odd integer and let $f$ be a Morse 
function on $M$ such that $\partial M=N_{a+}$. Then the above relation reads   
\[\sum_{k=0}^n(-1)^{n-k}c_k=\chi(M)+\chi(\partial M).\]
Now consider the Morse function $-f$ for which $\partial M=N_{a-}$, we get 
the following relation 
\[-\sum_{k=0}^n(-1)^{n-k}c_k=\chi(M).\]
Comparing these two relations one gets the following relation proving that the 
parity of the Euler character is a cobordism invariant.
\[\chi(\partial M)=-2\chi(M).\]
For example the Euler characteristic of a manifold which is the boundary 
of a contractible odd dimensional manifolds, e.g. $S^{2p}$ equals $-2$. 
A similar argument gives the relation 
$\chi(\partial M)=2\chi(M,\partial M)$.
Of course these relations can be obtained by homotop theory and Poincare 
duality. is a result of Poincare duality. 

{\bf Comment: } After diffusing this paper, the author has been informed by Maxim Braverman 
about the paper \cite{BrSi-Ki}, where a related and much more general result is proved also 
by Witten deformation method. In fact the Morse inequalities of this paper can be 
deduced from the inequalities proved in the work of Braverman and Silantyev. 
Nevertheless this paper provides 
a direct and relatively simple analytic proof for a special case, which
may be still interesting. The author would like to thank M. Braverman for informing me about 
their work and for other useful comments.

\section{Analytical properties of deformed Laplacian with 
boundary conditions}
To give a proof for theorem \ref{morsein}, we begin by studying the Laplacian 
operator, and its deformation, on the de Rham complex endowed with  
boundary conditions.  
Let $\omega$ be a differential 
$k$-form on $M$. In collar neighborhood $U$ of $\partial M$ it takes 
the following form 
\[\omega_{|U}=\omega_1(u,y)+du\wedge\omega_2(u,y)\] where 
$\omega_i$'s are $u$-depending differential forms on $\partial M$. 
Differential form $\omega$ satisfies the 
relative boundary condition $B_r$ if 
\begin{equation}\label{relboun}
\omega_1(0,y)=0\hspace{5mm}\text{ and }\hspace{5mm}
\frac{\partial \omega_2}{\partial u}(0,y)=0.
\end{equation}
Differential form $\omega$ satisfies the absolute boundary condition $B_a$ if 
\begin{equation}\label{absboun}
\frac{\partial\omega_1}{\partial u}(0,y)=0\hspace{5mm}\text{ and }\hspace{5mm}
\omega_2(0,y)=0.
\end{equation}
Clearly $\omega$ satisfies $B_r$ if and only if $*\omega$ satisfies $B_a$, 
where $*$ is a locally defined Hodge star operator with respect to a 
local orientation. 
>From now on we impose the relative boundary condition on $N_r$ and the 
absolute boundary condition on $N_a$ and denote this setting of 
boundary condition by $B$. By performing the completion of the set of 
smooth differential forms satisfying the boundary conditions $B$ with respect 
to appropriate Sobolev norm, we obtain the Sobolev spaces 
$W^l(M,\Lambda^*T*M;B)$. The inclusion $
W^l(M,\Lambda^*T*M;B)\hookrightarrow W^l(M,\Lambda^*T*M)$ is an isometry 
with closed image. So all classical theorems in the theory of the Sobolev 
spaces, $W^l(M,\Lambda^*T*M)$ e.g. 
Sobolev embedding theorem, Rellich's theorem and the elliptic estimate 
hold also in this context. 

Let $\eta$ be an another differential form on $M$ taking the 
form $\eta=\eta_1+du\wedge\eta_2$ in collar neighborhood of $\partial M$. 
The Green formula for $D:=d+\delta$(c.f \cite[]{}) takes the following form 
in our context  
\begin{equation}\label{gen}
\langle D\,\omega,\eta \rangle-\langle\omega,D\,\eta \rangle
=\int_{\partial M}\langle\omega_1,\eta_2 \rangle
-\int_{\partial M}\langle\omega_2,\eta_1 \rangle.
\end{equation}
So the operator $D$ is formally self adjoint
provided $\omega$ and $\eta$ satisfy both the relative or the absolute 
boundary conditions. Consequently $D$ is a formally self adjoint operator 
on $\Omega^*(M,B)$ and Laplacian operator $\triangle:=D^2$ is a 
formally positive second order elliptic operator on $\Omega^*(M,B)$ and 
\[\ker\triangle=\ker D.\]
The heat kernel of $e^{-t\triangle}$ may be constructed by means of 
heat kernels of boundary problems on half cylinder 
$\R^{\geq0}\times \partial M$ and heat kernel of Laplacian on closed 
manifold $M\sqcup_{\partial M}M$ (see \cite[page 55]{APS1}).
Since the function $f$ 
satisfies the Von Neumann condition \eqref{vonboun}, the multiplication by 
$e^{\pm sf}$ preserves the boundary condition $B$. Here, and so on, $s$ is a 
non negative real number. Let $d_s:=e^{-sf}d\,e^{sf}$ and  
$\delta_s=e^{sf}\delta e^{-sf}$ and put 
\[D_s:=d_s+\delta_s\hspace{5mm}\text{ and }\hspace{5mm}\triangle_s:=D_s^2.\]
As in above $D_s$ is formally self adjoint, so  the Witten Laplacian 
$\triangle_s$ is a formally positive elliptic differential form and 
\begin{equation}\label{kereq}
\ker\triangle_s=\ker D_s=\ker d_s\cap\ker\delta_s.
\end{equation}
The Hessian of the function $f$ is the following 2-tensor 
\[H(X,Y):=X.(Y.f)-(\nabla_XY).f\]
Let $\{e_i\}$ be a local orthonormal basis for $TM$ and $\alpha\in\Lambda^kT^*M$. 
Put $L_{e_i}(\alpha):=e^i.\alpha$ 
and $R_{e_i}(\alpha):=(-1)^k\alpha.e^i$, where dot ``.'' denotes the 
Clifford multiplication. The following relation defines a smooth 
section of $\End(\Lambda^*T^*M)$ which is independent of the orthonormal 
basis \[\mathsf H=\sum_{i,j}H(e_i,e_j)L_{e_i}R_{e_j}.\]
With these notation we have the following relation \cite[proposition 9.17]{Roe} 
\begin{equation}\label{moz}
\triangle_s^2=\triangle^2+s\mathsf H+s^2|df|^2.
\end{equation}
Notice that the operator $\mathsf H$ preserves clearly the boundary 
condition $B$. 
Moreover, due to the relation \eqref{vonboun} the multiplication operator $|df|$ 
preserve this boundary condition too.  
Using \eqref{moz} and the 
Duhamel formula (see, e.g. \cite{BeGeVe}), one can construct the heat operator 
$e^{-t\triangle_s}$ by means of $e^{-t\triangle}$. 
The heat operator is a smoothing operator on $L^2(M,\Lambda^*T^*M)$, 
so the Sobolev embedding theorem implies that it is 
a compact 
self adjoint operator on this Hilbert space. So, for each $0\leq k\leq n$, 
the Hilbert space 
$L^2(M,\wedge^kT^*M)$ may 
be decomposed to orthogonal sum of eigenspaces of $\triangle^k_s$. 
\begin{equation}\label{decsp}
L^2(M,\wedge^kT^*M)=\overline{\oplus_{\lambda_s\geq0}E^k_{\lambda_s}}.
\end{equation}
The following sequence is exact for $0\neq\lambda_s\in spec(\triangle_s)$ 
\begin{equation}\label{exeig}
0\to E^0_{\lambda_s}\stackrel{d_s}{\rightarrow}E^1_{\lambda_s}
\stackrel{d_s}{\rightarrow}
\dots\stackrel{d_s}{\rightarrow} E^{n-1}_{\lambda_s}
\stackrel{d_s}{\rightarrow} E^n_{\lambda_s}\to0.
\end{equation}
Notice that as a part of assertion, the exterior differential preserves 
the boundary conditions when it is restricted to eigenspaces of $\triangle_s$. 
We refere to \cite[page 49]{Gil} for a proof of this assertion. 
Since $d_s$ preserves $\ker \triangle$ the above assertion implies that 
the cohomology of the deformed de Rham complex
\[0\to\Omega^0(M,B)\stackrel{d_s}{\rightarrow}\Omega^1(M,B)
\stackrel{d_s}{\rightarrow}\dots\stackrel{d_s}{\rightarrow}\Omega^n(M,B)\to0,\] 
is isomorphic to $\ker\triangle_s$ and the isomorphism is induced from equality 
\eqref{kereq}. Moreover the multiplication by $e^{sf}$ provides an isomorphism 
between this complex and the ordinary de Rham complex corresponding to $s=0$. 
The cohomology of de Rham complex is the relative cohomology group 
$H^*_{dr}(M.N_r)$. So, as an immediate consequence of the above discussion 
we get the following relation
\begin{equation}\label{eqbetker}
\beta_k=\dim\ker\triangle_s^k.
\end{equation}
The basic analytical tool in the next section is the 
finite propagation speed property of wave operator that we shall to explain. 
The spectral resolution 
\eqref{decsp} can be used to define the wave operator $e^{itD_s}$. 
The restriction of this operator to $E^p_{\lambda_s}$ is the multiplication 
by $e^{it\lambda_s}$. With this construction of wave operator, it is clear that 
if $\omega$ is a smooth differential form with 
boundary condition $B$, then $\omega_t=e^{itD_s}\omega$ satisfies the 
boundary condition $B$ too, moreover it satisfies the wave equation 
\[(\frac{d}{dt}-iD_s)\omega_t=0.\] 
Let $\omega_t$ be a smooth solution of wave equation with $\omega_0=0$. 
The following relation shows that $\|\omega\|^2$ is independent of $t$
\[\frac{d}{dt}\|\omega\|^2=\langle iD_s\omega_t,\omega_t\rangle_2
+\langle\omega_t,iD_s\omega_t\rangle_2=0.\]
So $\omega_t=0$. Therefore the smooth solutions of wave equation are 
determined uniquely by their initial value.
\begin{lem}\label{finpro}
Let $\omega$ be a  compactly supported smooth differential p-form 
satisfying the boundary condition $B$. The support of 
$\omega_t=e^{itD_s}\omega$ is inside 
the distance $|t|$ of $supp(\omega)$. 
\end{lem}
\begin{pf}
>From the construction of wave operator one has  
$e^{tiD_s+t'iD_s}=e^{tiD_s}e^{t'iD_s}$ and $(e^{tiD_s})^*=e^{-tiD_s}$. 
So it suffices to prove the lemma for small positive values of $t$. 
For $x\in M$, Let $B(x,r)$ be a small geodesic ball around $x$ of radius $r$ and 
$S(x,r)$ be the unite sphere of radius $r$ which may intersect $\partial M$.  
As in the proof of \cite[proposition 7.20]{Roe}, the claim of the lemma 
follows from the energy estimate which asserts that for each $x\in M$
\[\int_{B(x,r-t)}|\omega_t|^2\]
is a decreasing function of $t$,where $d\mu$ is the Riemannian measure. 
In fact 
\[\frac{d}{dt}\int_{B(x,r-t)}|\omega_t|^2\,d\mu
=\int_{B(x,r-t)}\left(\langle iD_s\omega_t,\omega_t\rangle
+\langle\omega_t,iD_s\omega_t\rangle\right)d\mu-\int_{S(x,r-t)}|
\omega_t|^2\,d\sigma.\]
Here $d\sigma$ is the Riemannian measure induced on $S(x,r-t)$ as the 
boundary of $B(x,r-t)$.
Using the Green formula, the first integral on the right side of above 
relation equals 
\[\int_{\partial B(x,r-t)}\langle\rm{cl}(n).\omega_t,\omega_t\rangle d\sigma
=\int_{S(x,r-t)}\langle\rm{cl}(n).\omega_t,\omega_t\rangle d\sigma
+\int_{B(x,r-t)\cap\partial M}\langle\rm{cl}(n).\omega_t,
\omega_t\rangle d\sigma.\]
It follows from relation \eqref{gen} that the points on $\partial M$ 
have no contribution in the above integral, due to the boundary conditions. 
In other hand 
the Clifford action of the normal vector $n$ to $S(x,r-t)$ is 
a pointwise isometry, so the Cauchy inequality implies that 
\[|\int_{S(x,r-t)}\langle\rm{cl(n)}.\omega_t,\omega_t\rangle d\sigma|
\leq \int_{S(x,r-t)}|\omega_t|^2\,d\sigma.\]
This shows $\frac{d}{dt}\int_{B(x,r-t)}|\omega_t|^2\,d\mu\leq0$ and completes 
the proof.
\end{pf}

\section{derivation of the Morse inequalities}
Let $\phi$ be a non negative rapidly decreasing 
function on $\R^{\geq0}$ satisfying $\phi(0)=1$. 
The operator $\phi(\triangle_s^k)$ being a smoothing operator,
is of trace class, so we can define 
$\nu_k(s):=\rm{Tr}\,\phi(\triangle_s^k)$. 
The following proposition is a generalized version of the proposition $14.3$ of  
\cite{Roe} in the presence of the boundary. We prove this theorem by using 
the spectral resolution 
provided by eigenvectors of $\triangle_s$.
\begin{prop}\label{anamor}
\label{anamorin} The following inequalities hold for 
$0\leq k\leq n$
\begin{equation}
\nu_k(s)-\nu_{k-1}(s)+\dots+(-1)^k\nu_0(s)\geq\beta_k-\beta_{k-1}+
\dots+(-1)^k\beta_0.
\end{equation}
and equality holds for $k=n$. 
\end{prop}
\begin{pf}
Consider the following exact sequence coming from \eqref{exeig}
\begin{equation}
0\to E^0_{\lambda_s}\stackrel{d_s}{\rightarrow}E^1_{\lambda_s}
\stackrel{d_s}{\rightarrow}
\dots\stackrel{d_s}{\rightarrow} E^{k}_{\lambda_s}
\stackrel{d_s}{\rightarrow}\rm{Im}\,d_s^k\to0.
\end{equation}
Obviously $\phi(\triangle_s)$ restricts to a linear operator on $\rm{Im}\,d_s^k$. 
The trace of this restriction, denoted by $r(\phi,\lambda_s,k)$, is non negative 
for $0\leq k\leq n$ and is zero for $k=n$. 
Since $\triangle_{s|E_{\lambda_s}}=\lambda_s\,\id$, the trace of the
restriction of $\phi(\triangle_s)$ to $E_{\lambda_s}^j$ equals 
$\phi(\lambda_s).\dim E_{\lambda_s}^j$. So 
\[r(\phi,\lambda_s,k)-\rm{tr}\,
\phi(\triangle^k_{\lambda_s|E^k_{\lambda_s}})+\rm{tr}\,
\phi(\triangle^{k-1}_{\lambda_s|E^{k-1}_{\lambda_s}})
-\dots+(-1)^k\rm{tr}\,\phi(\triangle^0_{\lambda_s|E^0_{\lambda_s}})=0.\]
By summation over all $\lambda_s\neq0$ and using the fact that $\phi(0)=1$ 
and \eqref{eqbetker} we get 
\[\sum_{\lambda_s\neq0}r(\phi,\nabla_s,k)
-(\mu_k(s)-\beta_k)+(\mu_{k-1}(s)-\beta_{k-1})-\dots
+(-1)^k(\mu_0(s)-\beta_0)=0.\]
These inequalities prove the assertion of the proposition because for 
$0\leq k\leq n$ one has 
$r(\phi,\lambda_s,k)\geq0$ and equality holds for $k=n$. 
\end{pf}
With above proposition,  
to prove the Morse inequalities of theorem \ref{morsein} 
we have to study the traces 
$\nu_k(s)=\rm{Tr}\,\phi(\triangle_s^k)$ when $s$ goes to infinity. 
For a small positive number $\rho<1/4$, let $M_\rho$ denote 
the disjoint union of $\rho$-neighborhoods of critical points and the 
$\rho$-collar neighborhood $(\rho,0]\times\partial M$ of the boundary. 
The positive number $\rho$ is so small that each connected component of 
$M_{4\rho}$ contains only one critical point or one connected component 
of boundary $\partial M$. We recall that $\phi(\triangle_s^k)$ is a 
smoothing operator with smooth kernel 
$K(s,x,x')\in\wedge^kT_xM\otimes\wedge^kT_{x'}M$. 
 
\begin{prop}\label{noncritcon}
Let $\phi$ be a rapidly decreasing even 
function such that the fourier transform of function $\psi$ defined 
by $\psi(t):=\phi(t^2)$ is supported in $(-\rho,\rho)$.
When $s$ goes toward infinity, the kernel $K(s,x,x)$ 
of $\phi(\triangle_s^k)$ goes 
uniformely to $0$ for $x$ in  $M\verb+\+M_{2\rho}$. 
\end{prop}
\begin{pf}
In the complement of $M_{4\rho}$ one has $|df|\geq c$ for some positive 
constant  $c$. So, using the relation \eqref{moz}, one get 
the following estimate  
\begin{equation}\label{princ}
 \langle D_s^2\omega,\omega\rangle_2\geq\frac{1}{2}cs^2\|\omega\|_2^2
\hspace{5mm}\text{ if }\hspace{3mm}supp(\omega)\subset M_{4\rho}.
\end{equation}
Using this inequality, the finite propagation speed property of lemma 
\ref{finpro} and classical theorem of Sobolev spaces $W^k(M,\Lambda^**M;B)$, 
the proof of the proposition 14.6 in \cite{Roe} goes over verbatim and prove 
this proposition.
\end{pf}
Let $\beta$ be a smooth function on $M$ which is supported in $M_{3\rho}$ 
and is equal to $1$ on $M_{2\rho}$. The above lemma shows that 
\begin{equation}\label{loccon}
\lim_{s\to\infty}\rm{Tr}\,\phi(\triangle_s^k)
=\lim_{s\to\infty}\rm{Tr}\,\,\bar\beta\phi(\triangle_s^k),
\end{equation}
where $\bar\beta$ is the pointwise multiplication by $\beta$. 
So, the next step is to study 
the asymptotic behavior of $\rm{Tr}\,\,\bar\beta\phi(\triangle_s^k)$ 
when $s$ goes 
to infinity. From now on we consider $M_{4\rho}$ either with coordinates given by 
Morse lemma around the non degenerate critical points, or  
by collar coordinates $(u,y)$ in $(4\rho,0]\times\partial M$. 
We assume that the Riemannian structure is euclidian around 
critical points and is of product form $du^2+g_0$ in the 
collar neighborhood. Using these coordinates   
we can identify the differential form supported in $M_{4\rho}$ 
with differential forms supported either in the $4\rho$-neighborhood of origin 
in $\R^n$ or in product space $(4\rho,0]\times\partial M$. 
In the neighborhood of a critical point $x_0\in M$ with Morse index $r$, 
the function $f$ takes the following form 
\begin{equation}\label{locmor}
f(x_1,\dots,x_n)=f(x_0)-\frac{1}{2}x_1^2-\dots-\frac{1}{2}x_r^2+
\dots +\frac{1}{2}x_n^2.
\end{equation}
So the deformed Laplacian $\triangle_s^k$ given by relation \eqref{moz} 
coincides with the following operator acting on $\Omega^k(\R^n)$
\begin{equation}\label{harmon0}
L_s^k=-\sum_{j=1}^n\frac{\partial^2}{\partial^2x_j}
+s^2x_j^2\,+s\epsilon_jZ_j.
\end{equation}
Here $\epsilon_j=\pm1$ is the sign of the coefficient of $x_j^2$ in the expression \eqref{locmor}. 
Moreover $Z_j:=[dx_j\wedge.,dx_j\llcorner.]$. In the collar neighborhood 
$[0,1)\times N_+$ we have $f(u,y)=\frac{1}{2}u^2$. So the expression \eqref{moz} 
for $\triangle_s^k$ gives the following operator acting on 
$\oplus_{\epsilon=0,1}\Omega^\epsilon([0,1),B)
\otimes\Omega^{k-\epsilon}(N_+)$
\begin{equation}\label{euclapboun+0}
L_s^k=(A_s^0+\triangle^k_{N_+})\oplus(A_s^1+\triangle^{k-1}_{N_+}),
\end{equation}
where 
\[A_s^\epsilon=-\frac{\partial^2}{\partial u^2}+s^2u^2+(-1)^{\epsilon+1} s.\] 
In the collar neighborhood $[0,1)\times N_-$ we have $f(u,y)=-\frac{1}{2}u^2$. 
So the expression \eqref{moz} 
for $\triangle_s^k$ gives the following operator acting on  
$\oplus_{\epsilon=0,1}\Omega^\epsilon([0,1),B)
\otimes\Omega^{k-\epsilon}(N_-)$ 
\begin{equation}\label{euclapboun-0}
L_s^k=(A_s^0+\triangle^k_{N_-})\oplus(A_s^1+\triangle^{k-1}_{N_-}),
\end{equation}
where
\[A_s^\epsilon=-\frac{\partial^2}{\partial u^2}+s^2u^2+(-1)^{\epsilon} s,\]
In above discussion $\triangle^{k-1}_{N_+}$ and $\triangle^{k-1}_{N_-}$ 
are, respectively, the Laplacian operators on $N_+$ and $N_-$.
\begin{prop}\label{proploc}
Let $L_s^k$ denote any one of operators given by 
\eqref{harmon0},\eqref{euclapboun+0} or \eqref{euclapboun-0} and let 
$\bar\beta$ be the 
corresponding operator defined in above. The following equality holds 
\begin{equation*}
\rm{Tr}\,\,\bar\beta\phi(\triangle_s^k))=\rm{Tr}\,\,\bar\beta\phi(L_s^k).
\end{equation*}
Here $\phi(\triangle_s^k)$ and $\phi(L_s^k)$ are bounded operator on 
different $L^2$-Hilbert spaces.
\end{prop}
\begin{pf}
Clearly $\triangle_s=D_s^2$ and $L_s=A_s^2$ where $A_s$ is a 
differential operator 
which is equal to $D_s$ in $M_{4\rho}$. Define the function $\psi$ 
by $psi(t):=\phi(t^2)$. 
To prove the proposition it suffices to prove the following equality, 
provided that $\omega$ 
is supported in $M_{3\rho}$,
\[\psi(D_s)\,\omega=\psi(A_s)\,\omega.\]
For this purpose consider the following partial differential equation for 
 t-depending 
differential forms which vanish in $M\verb+\+M_{3\rho}$
\[\frac{\partial^2\omega_t}{\partial^2t}+D^2_s\,\omega_t=0.\]
It is easy to verify that the {\it energy} of a smooth solution $\omega_t$,
\[\|\frac{\partial \omega_t}{\partial t}\|^2+\|D_s\omega_t\|^2,\]
is independent of $t$. So, the solution is uniquely determined by 
initial conditions 
$\omega_{t|t=0}$ and $(\frac{\partial\omega_t}{\partial t})_{|t=0}$. 
Given a smooth differential 
form $\omega$ with support in $M_{2\rho}$, the relations 
$\omega_t:=\cos(tD_s)\omega=\frac{1}{2}(e^{itD_s}+e^{-itD_s})$ and 
 $\omega'_t:=\cos(tA_s)$ are 
clearly  smooth solutions of above differential equation with the same 
initial conditions. 
Therefor the unite propagation speed property of lemma \ref{finpro} implies 
the equality 
$\omega_t=\omega'_t$ for $|t|<\rho$. Using the fact that $\hat\psi$ is an even 
 function which is 
supported in $[-\rho,\rho]$, we have
\begin{align*}
\psi(D_s)\omega&=\frac{1}{2\pi}\int_{-\rho}^\rho\hat\psi(t) e^{itD_s}\omega\,dt\\
&=\frac{1}{\pi}\int_0^\rho\hat\psi(t)\cos(tD_s)\omega\,dt\\
&=\frac{1}{\pi}\int_0^\rho\hat\psi(t)\cos(tL_s)\omega\,dt\\
&\vdots\\
&=\psi(A_s)\omega.
\end{align*}
This completes the proof of the proposition.
\end{pf}
In view of the above proposition, for computing $\nu_k(s)$ when $s$ goes 
toward $\infty$, it suffices to compute $Tr(\bar\beta\phi(L_s^k))$, at 
$s=\infty$, where $L_s^k$ is given by relations \eqref{harmon0}, 
\eqref{euclapboun+0} and \eqref{euclapboun-0}.

Concerning the operator $L_s^k$ given by \eqref{harmon0}, the following relation 
holds, cf. \cite[lemma 14.11]{Roe}
\begin{equation}\label{crit}
\lim_{s\to\infty}\rm{Tr}(\bar\beta\,\phi(L_s^k))
=\lbrace\begin{array}{cc}
0&k\neq r\\
1&k=r
\end{array}
\end{equation}
Now we are going to compute 
$\lim_{s\to0}\rm{Tr}\,(\bar\beta\phi(L_s^k))$ for $L_s^k$ of relations 
\eqref{euclapboun+0} and \eqref{euclapboun-0}. 
For this purpose, we summarize some basic properties of 
the {\it one dimensional harmonic oscillator operator}(see \cite[chapter 9]{Roe}) 
\[-\frac{\partial^2}{\partial^2u}+s^2u^2\colon\mathcal{S}(\R)
\to \mathcal{S}(\R),\]
where $\mathcal{S}(\R)$ denotes the space of smooth rapidly decreasing 
functions on $\R$. 
The eigenvalues of this operator are $(2p+1)s$ for $p=0,1,2,\dots$ and 
each eigenvalue has multiplicity one. 
The coresponding eigenfunctions have the following form 
\[\vartheta_p(s,u)=(2^ps^{p+1}p!)^{-\frac{1}{2}}(-\frac{d}{du}+su)^p
e^{-\frac{su^2}{2}}.\]
In particular these eigenfunctions satisfy the Dirichlet boundary 
condition when $p$ is odd and the Von Neumann condition when $p$ is even. 
Another basic property is related to $\vartheta_0$. If $\beta$ is a 
rapidly decreasing continuous function defined around $0\in\R$ then 
\begin{equation}\label{ator}
\lim_{s\to\infty}\langle(\beta(.)\vartheta_0(s,.)),\vartheta_0(s,.)\rangle
=\beta(0). 
\end{equation}
Obviously the harmonic oscillator operator may be considered as an operator on 
the Schwartz space 
$\mathcal{S}(\R)\,du$ of differential $1$-forms. The eigenvalues of this 
operator are the same and 
the eigenvectors are $\vartheta_p(s,u)\,du$.
\newline
In what follows we will need to consider a sequence $\phi_m;~~m\in\N$ 
of rapidley decreasing non negative even functions on $\R$ satisfying  
$\phi_m(0)=1$ such that the Fourier transforms $\hat\psi$ is supported 
in $(-\rho,\rho)$ where $\psi_m(t)=\phi(t^2)$. Moreover the sequence $\phi_m$ 
converges uniformly to zero, out of the compact neighborhoods of $0\in\R$.  

\subsubsection{{\bf Cylindrical operators arising from $N_+$.}}
On the product space $\R^{\geq0}\times N_+$ we consider the deformed 
Laplacian operator $L_s^k=A_s^\epsilon+\triangle_{N_+}$ given by 
\eqref{euclapboun+0} and acting on 
$\Omega^\epsilon([0,1),B_r)\otimes\Omega^{k-\epsilon}(\partial M)$. Here 
\begin{equation}
A_s^\epsilon=-\frac{\partial^2}{\partial u^2}+s^2u^2+(-1)^{\epsilon+1} s,
\end{equation}
Let $\{\psi_{\lambda_k}^k\}_{\lambda_k}$ be a spectral resolution for $\triangle_s^k$.
So the eigenvectors of the 
restriction of $L_s^k$ to $\R^{\geq0}\times N_{r+}$ with respect to boundary 
condition \eqref{relboun} are 
\[\vartheta_{2l+1}\otimes\psi^k_{\lambda_k}\hspace{5mm};\hspace{5mm}
\vartheta_{2l}\,du\otimes
\psi^{k-1}_{\lambda_{k-1}},\hspace{1cm}l=0,1,2,\dots.\] 
so 
\begin{align*}
\rm{Tr}\,\,\bar\beta\phi_m(L^k_{s|\R^{\geq0}\times N_{r+}})
&=\sum_{l,\,\lambda_k}\phi_m(s+4ls+\lambda_k)\langle\beta(.)
\vartheta_{2l+1}(s,.),\vartheta_{2l+1}(s,.)\rangle\\
&+\sum_{l,\lambda_{k-1}}\phi_m(2s+4ls+\lambda_{k-1})\langle\beta(.)
\vartheta_{2l}(s,.),\vartheta_{2l}(s,.)\rangle.
\end{align*}
The function $\phi_m$ is rapidly decreasing at infinity and the arguments of 
$\phi_m$ appearing in above formula go to infinity when $s$ do, so 
\begin{equation}\label{rel+}
\lim_{s\to\infty}\rm{Tr}\,\,\bar\beta\phi_m(L^k_{s|\R^{\geq0}\times N_{r+}})=0.
\end{equation}
Now we consider the case of $L^k_{s|\R^{\geq0}\times N_{a+}}$, given 
by \eqref{euclapboun-0}, with respect 
to absolute boundary conditions \eqref{absboun}. According to above 
discussion, the eigenvectors of this operators are 
\[\vartheta_{2l}\otimes\psi^k_{\lambda_k}\hspace{5mm};\hspace{5mm}
\vartheta_{2l+1}\,du\otimes
\psi^{k-1}_{\lambda_{k-1}},\hspace{1cm}l=0,1,2,\dots.\]
So we have 
\begin{align*}
\rm{Tr}\,\,\bar\beta\phi_m(L^k_{s|\R^{\geq0}\times N_{r+}})
&=\sum_{l,\,\lambda_k}\phi_m(4ls+\lambda_k)\langle\beta(.)
\vartheta_{2l}(s,.),\vartheta_{2l}(s,.)\rangle\\
&+\sum_{l,\lambda_{k-1}}\phi_m(3s+4ls+\lambda_{k-1})\langle\beta(.)
\vartheta_{2l+1}(s,.),\vartheta_{2l+1}(s,.)\rangle.
\end{align*}
The argument of $\phi_m$ is second summation appearing in above equalities go  
all toward infinity when $s$ do, so this sum has no  contribution 
when $s$ goes to infinity. The arguments of $\phi_m$ in the first 
summation, except those corresponding to $l=\lambda_k=0$, have no 
contribution because $\phi_m$'s are rapidly decreasing at 
infinity and because the sequence $\{\phi_m\}_m$ converges uniformly to zero 
out of compact neighborhoods of $0\in\R$. Therefore 
\begin{equation*}
\lim_{s,m\to\infty}\rm{Tr}\,\,\bar\beta\phi_m(L^k_{s|\R^{\geq0}\times N_{a+}})
=\lim_{s,m\to\infty}\sum_{\lambda_k=0}\phi_m(0)\langle\beta(.)
\vartheta_{0}(s,.),\vartheta_{0}(s,.)\rangle.
\end{equation*}
Now using the relation \eqref{ator} and using the Hodge theory on $N_{a+}$ we 
conclude 
\begin{equation}\label{abs+}
\lim_{s\to\infty}\rm{Tr}\,\,\bar\beta\,\phi_m(L^k_{s|\R^{\geq0}
\times N_{a+}})=\eta_k.
\end{equation}

\subsubsection{{\bf Cylindrical operators arising from $N_-$.}}
On the product space $\R^{\geq0}\times N_-$ we consider the deformed 
Laplacian operator $L_s=A_s^\epsilon+\triangle_{\partial M}^j$ given 
by \eqref{euclapboun-0} and acting on 
$\Omega^\epsilon([0,1),B_r)\otimes\Omega^j(\partial M)$. Here  
\begin{equation}\label{euclapboun-}
A_s^\epsilon=-\frac{\partial^2}{\partial u^2}+s^2u^2+(-1)^{\epsilon} s.
\end{equation}
The eigenvectors of $L_s^k$ with respect to boundary condition 
\eqref{relboun} are 
\[\vartheta_{2l+1}\otimes\psi^k_{\lambda_k}\hspace{5mm};\hspace{5mm}
\vartheta_{2l}\,du\otimes
\psi^{k-1}_{\lambda_{k-1}},\hspace{1cm}l=0,1,2,\dots.\] 
therefore
\begin{align*}
\rm{Tr}\,\,\bar\beta\phi_m(L^k_{s|\R^{\geq0}\times N_{r+}})
&=\sum_{l,\,\lambda_k}\phi_m(3s+4ls+\lambda_k)\langle\beta(.)
\vartheta_{2l+1}(s,.),\vartheta_{2l+1}(s,.)\rangle\\
&+\sum_{l,\lambda_{k-1}}\phi_m(4ls+\lambda_{k-1})\langle\beta(.)
\vartheta_{2l}(s,.),\vartheta_{2l}(s,.)\rangle.
\end{align*}
All arguments of $\phi_m$ in the first summation go to infinity when $s$ do, 
so this summation has no non-vanishing  contribution at $s=\infty$. 
The sequence $\phi_m$ converges uniformly to zero out of compact 
neighborhoods of $0\in\R$. So in the second summation, only the terms 
corresponding to $l=\lambda_{k-1}=0$ may have non vanishing contribution 
when $s$ and $m$ go to infinity. Therefore using relation \eqref{ator} 
\begin{align}
\lim_{s,m\to\infty}\rm{Tr}\,\,\bar\beta\phi_m(L^k_{s|\R^{\geq0}\times N_{r-}})
&=\lim_{s,m\to\infty}\sum_{\lambda_{k-1}=0}\phi_m(0)\langle\beta(.)
\vartheta_{0}(s,.),\vartheta_{0}(s,.)\rangle\notag\\
&=\gamma_{k-1}\label{rel-}.
\end{align}
A similar discussion gives the following relation 
\begin{equation}\label{abs-}
\lim_{m,s\to\infty}\rm{Tr}\,\,\bar\beta\phi_m(L^k_{s|\R^{\geq0}\times N_{a-}}))=0.
\end{equation}
Now we are ready to give the proof of theorem \ref{morsein}

{\bf proof of theorem \ref{morsein}}:. 
Put $\mu_k(m,s):=\rm{Tr}\,\phi_m(L_s^k)$. Let $\beta$ be a continuous 
function supported in small neighborhoods of critical points of $f$ and in 
collar neighborhood $(1,0]\times\partial M$ of boundary. The value of 
this function on critical points and on boundary is assumed to be $1$. 
Moreoner in collar neighborhood $(1,0]\times\partial M$ it is assumed to 
be a function of $u$. We denote by 
$\bar\beta$ the pointwise multiplication of differential forms by $\beta$. 
>From relation 
\eqref{loccon} we get 
\[\lim_{s\to\infty}\nu_k(m,s)
=\lim_{s\to\infty}\rm{Tr}\,\,\bar\beta\phi_m(\triangle_s^k).\]
Proposition \ref{proploc} implies 
\[\rm{Tr}\,\bar\beta\phi_m(\triangle_s^k)=\sum\rm{Tr}
\,\bar\beta\phi_m(L_s^k),\]
where the sum is taken over connected component of $M_{4\rho}$.
Now relations \eqref{crit},\eqref{rel+},\eqref{abs+},\eqref{rel-} 
and \eqref{abs-} together imply 
\[\lim_{s,m\to\infty}\sum\rm{Tr}\,\,\bar\beta\phi_m(L_s^k)
=c_k+\eta_k+\gamma_{k-1}.\]
This relation with analytic Morse inequalities given in proposition  
\ref{anamorin} give the desired results of the theorem.
\begin{rem}
At the beginning of the paper we assumed that the Morse function 
$f\colon M\to \R$ takes the form 
$f(u,y)=\pm\frac{1}{2}u^2$ in the collar neighborhood $[0,1)\times\partial M$. 
This condition may 
be weakened and we could consider smooth functions which are Morse both in 
the interior and on 
the boundary of $M$. The restriction of such a function, $f\colon M\to\R$, 
to the boundary is a 
Morse function too and the inductive proof of the Morse lemma 
(see \cite{Mi1}) shows that the
Morse functions $f$ takes the following form with respect to a local coordinates 
$u,y_1,\dots,y_{n-1}$ 
\[f(u,y)=f(x_0)\pm\frac{1}{2}u^2\pm\frac{1}{2}y_1^2\dots\pm\frac{1}{2}y_{n-1}^2.\] 
Here $(y_1\dots,y_{n-1})$ is a local coordinates system for $\partial M$ around 
the critical 
point $x_0\in\partial M$. For a connected component $N$ of the boundary 
$\partial M$, 
let $c_{k,k'}(N)$ denote the number of those critical points $x_0\in N$ of $f$  such 
that $\ind_{x_0}f=k$ and $\ind_{x_0}(f_{N}=k'$. As before let $c_k$ denote 
the number of 
critical points of $f$ with Morse index $k$ which are in the interior of $M$. 
Put $\mu_k:=c_k+c_{k,k}(N_a)+c_{k,k-1}(N_r)$ and $\beta_k=\dim H_{dr}(M,N_r)$. 
The above discussion about cylindrical operators and boundary properties 
of eigenfunctions of 
the harmonic oscillator can be applied to the operator $L_s^k$ appeared in 
relation \eqref{crit} to 
deduce the following Morse inequalities
\begin{equation*}
\mu_k-\mu_{k-1}+\dots+(-1)^k\mu_0\geq\beta_k-\beta_{k-1}+
\dots+(-1)^k\beta_0.
\end{equation*} 
As before the equality holds for $k=n$. 
\end{rem}

\end{document}